\documentclass{birkjour}
\usepackage{amssymb}
\usepackage{amsmath}
\usepackage{amsfonts}
\theoremstyle{plain}
\newtheorem{thm}{Theorem}[section]

\newtheorem{lmm}[thm]{Lemma}
\newtheorem{cor}[thm]{Corollary}

\newtheorem{rmk}[thm]{Remark}

\newtheorem{prb}[thm]{Problem}
\theoremstyle{remark}

\def\II{\mathbb{I}}
\def\pmc#1{\setbox0=\hbox{#1}
    \kern-.1em\copy0\kern-\wd0
    \kern.1em\copy0\kern-\wd0}
\def\al{\alpha}
\def\be{\beta}

\def\si{\sigma}
\def\Si{\Sigma}

\def\om{\omega}

\def\ov{\overline}

\begin{document}

\title[On the singular homology]{On the singular homology of one class of simply-connected cell-like spaces}

\author[K. Eda]{Katsuya Eda}
\address{School of Science and Engineering,
Waseda University, Tokyo 169-8555, Japan}
\email{eda@logic.info.waseda.ac.jp}

\author[U. H. Karimov]{Umed H. Karimov}
\address{Institute of Mathematics,
Academy of Sciences of Tajikistan,
Ul. Ainy $299^A$, Dushanbe 734063, Tajikistan}
\email{umedkarimov@gmail.com}

\author[D. Repov\v{s}]{Du\v san Repov\v s}
\address{Faculty of Mathematics and Physics, and
Faculty of Education,
University of Ljubljana, P.O.Box 2964,
Ljubljana 1001, Slovenia}
\email{dusan.repovs@guest.arnes.si}

\date{\today}



\begin{abstract}

In our earlier papers we constructed examples of
2-dimensional nonaspherical
simply-connected cell-like Peano continua, called {\sl  Snake
space}. In the sequel we introduced the functor $SC(-,-)$ defined on
the category of all spaces with 
base points and continuous
mappings. For the circle $S^1$, the space $SC(S^1, \ast)$ is a
Snake
space. In the present paper we study the higher-dimensional
homology and homotopy properties of the spaces
$SC(Z, \ast)$ for any path-connected compact spaces $Z$.
\end{abstract}

\date{\today}
\maketitle

\section{Introduction}

It is well-known that there exist planar noncontractible continua
$X$ all homotopy groups $\pi_i(X),\ i \geq 1,$ of which are
trivial (e.g. the {\it Warsaw circle}). Every planar simply
connected Peano continum is a
contractible space, see e.g. \cite{KRRZ,N}. 
Noncontractible homology locally connected 
(HLC and therefore
Peano) continua, 
all homotopy groups of which are trivial,
were constructed
in \cite{KR:Proc2009}. All these
examples are infinite-dimensional. The following
problem remains
open \cite{EKR:nonaspherical}:

\begin{prb}\label{1}
Does there exist a finite-dimensional noncontractible Peano
continuum all homotopy groups of which are trivial?
\end{prb}

 We constructed
 in \cite{EKR:topsin} the functor
$SC(-,-)$,  defined on the category of all topological spaces with
base points.
Roughly speaking, for any space $Z$, one takes  the
infinite cylinder $Z\times [0,\infty )$ and attaches it to the
square $[0,1]\times[-1,1]\subset \mathbb{R}^{2}$, along the open
{\it Topologist sine curve}:
$$
\{ (x,y) \in \mathbb{R}^{2} \mid  y= \sin(1/x), \ 0<x\le 1\},$$ 
so
that its diameter tends to zero. The space
$SC(Z, \ast)$ is
called the {\it Snake cone} and when $Z$ is the circle $S^1$, the
space $SC(S^1, \ast)$ is called the
{\it Snake space}.

The Snake space was the first candidate for an example of a
simply
connected aspherical noncontractible Peano continuum. However,
 we have discovered, rather unexpectedly, that
 the group
 $\pi_2(SC(S^1))$ is nontrivial
\cite{EKR:nonaspherical}.

It is easy to see that the Snake space is a
cell-like  Peano
continuum (for the verification
of cell-likeness use e.g. \cite{MS}). We have already proved  the following:

\begin{thm}\label{thm:main1} \cite[Theorem 1.1]{EKR:topsin}
For every path-connected space $Z,$ the Snake cone $SC(Z)$ is
simply-connected.
\end{thm}

Our original proof in \cite{EKR:topsin} was quite long and
technical. We shall give a short proof of this result in Section 2
of the present paper.

We proved
in \cite{EKR:nonaspherical}  that whenever
$\pi _1(Z, z_{0})$
is nontrivial, the singular homology group
$H_2(SC(Z);\mathbb{Z})$ is nontrivial,
and since the spaces $SC(Z)$
are simply connected,
it follows
by the Hurewicz Theorem,
that
$\pi _2(SC(Z),z_{0})$ is isomorphic to
$H_2(SC(Z);\mathbb{Z}),$
and hence is also nontrivial. 
The converse was proved in
\cite{EKR:glasnik}, along the lines of the proof in
\cite{EKR:topsin}. In Section 3 of the present paper we shall give
a different and significantly shorter proof of the generalized
result for $(n-1)-$connected spaces, when $n\geq 2$.

\begin{thm}\label{thm:main2}$($\cite[Theorem 1.1]{EKR:glasnik} for $n=2$.$)$
Let  $Z$ be any $(n\! -\! 1)$-connected space, $n\ge 2$. Then
$H_n(SC(Z);\mathbb Z)$ and $\pi _n(SC(Z), z_{0})$ are trivial.
\end{thm}
Since by \cite[Theorem 3.1]{EKR:nonaspherical} the nontriviality of $\pi
_1(Z)$ implies that of $H_2(SC(Z))$, we obtain
the following:
\begin{cor}
For any  path-connected space $Z$ and any point
$z_{0}\in Z$, the following statements are equivalent:\\
$({\rm i})$ $\pi _1(Z, z_{0})$  is trivial;\\
$({\rm ii})$ $H_2(SC(Z); \mathbb{Z})$  is trivial; and\\
$({\rm iii})$ $\pi _2(SC(Z), z_{0})$  is trivial.
\end{cor}

Undefined notions are the usual ones and we refer the reader to
\cite{Spanier:algtop}.

\section{Proof of Theorem~\ref{thm:main1}}

We shall
follow the notations  for the Snake cone $SC(Z)$ and
the projection $p: SC(Z)\to \II^2$ as in \cite{EKR:topsin}.
The polygonal line $A_1B_1A_2B_2\cdots $ on $\II^2$,
together with the limit interval $AB$,
is the
{\it piecewise linear version} of the Topologist sine curve in Figure 1.

\begin{center}
\setlength{\unitlength}{0.7mm}
\begin{picture}(110,110)
\put(26,101){\line(-1,-6){10}}
\put(26,0){\line(0,1){101}}
\put(26,0){\line(1,4){25}}
\put(51,0){\line(0,1){101}}
\put(51,0){\line(1,2){50}}
\put(1,1){\line(1,0){100}}
\put(1,101){\line(1,0){100}}
\put(1,1){\line(0,1){100}}
\put(101,1){\line(0,1){100}}

\put(-3,-4){\shortstack{$A$}}
\put(-2,104){\shortstack{$B$}}
\put(-5,50){\shortstack{C}}
\put(102,104){\shortstack{$B_1$}}
\put(50,104){\shortstack{$B_2$}}
\put(25,104){\shortstack{$B_3$}}
\put(102,-4){\shortstack{$A_1$}}
\put(50,-4){\shortstack{$A_2$}}
\put(25,-4){\shortstack{$A_3$}}
\put(103,50){\shortstack{$C_1$}}
\put(78,50){\shortstack{$C_2$}}
\put(53,50){\shortstack{$C_3$}}
\put(40,50){\shortstack{$C_4$}}
\put(28,50){\shortstack{$C_5$}}

\put(100,0){\shortstack{$\bullet$}}
\put(100,50){\shortstack{$\bullet$}}
\put(75,50){\shortstack{$\bullet$}}
\put(50,50){\shortstack{$\bullet$}}
\put(37,50){\shortstack{$\bullet$}}
\put(50,0){\shortstack{$\bullet$}}
\put(25,0){\shortstack{$\bullet$}}
\put(25,50){\shortstack{$\bullet$}}
\put(25,100){\shortstack{$\bullet$}}
\put(50,100){\shortstack{$\bullet$}}
\put(100,100){\shortstack{$\bullet$}}
\put(0,100){\shortstack{$\bullet$}}
\put(0,50){\shortstack{$\bullet$}}
\put(0,0){\shortstack{$\bullet$}}
\end{picture}
\end{center}

\smallskip
\begin{center}
Figure 1
\end{center}
\smallskip

For a proof of Theorem~\ref{thm:main1}, we recall the
notion of the
free
$\si$-product of groups and a lemma from \cite{E:free}.
An element of $\pmc{$\times$}\, \, \; _{i\in I} ^\si G_i$ is expressed by a
word $W\in \mathcal{W}^\si (G_i:i\in I)$,
where $G_i\cap G_j =
\{ e\}$ for $i\neq j$, $W:\ov{W}\to \bigcup \{ G_i: i\in I\}$, $\ov{W}$
is a countable linearly ordered set, and $W^{-1}(G_i)$ is finite for
each $i\in I$
(cf. \cite{E:free}).
Let $h: \pmc{$\times$}\, \, \; _{i\in I} ^\si G_i \to
\pmc{$\times$}\, \, \; _{j\in J} ^\si H_j$ be a homomorphism.

For any $W\in \mathcal{W}^\si (G_i:i\in I)$,
express $h(W(\al ))$ by a reduced
word $V_\al \in \mathcal{W}^\si (H_j:j\in J)$.
Define $\ov{V}$ to be $\{ (\al ,\be ): \al \in \ov{W}, \be \in \ov{V_\al}\}$
with the lexicographical ordering and $V(\al ,\be) = V_\al (\be )$.
A homomorphism $h$ is said to be
{\it standard}, if $V$ as defined above, is
a word in $\mathcal{W}^\si (H_j:j\in J)$ and $h(W) = V$, for every $W
\in \mathcal{W}^\si (G_i:i\in I)$.
We have used the superscript ${}^\si$ in some cases, which means a
restriction to the countable case.
Hence, when an index set is countable,
the restriction is unnecessary and we drop the superscript ${}^\si$. 

Let $\{(X_i,x_i)\}_{i\in I}$ be pointed spaces. Let
$(\widetilde{\bigvee}_{i\in I}(X_i,x_i),x^*)$ be a
bouquet of
$\{( X_i,x_i)\}_{i\in I}$. The underlying set
$(\widetilde{\bigvee}_{i\in I}(X_i,x_i),x^*)$ is the quotient
space of a
discrete union of all $X_i$'s by the
identification of all
points $x_i $ with a singleton $x^*$ and the topology is defined
by specifying the neighborhood bases as follows (c.f.  \cite{BM}):

\begin{itemize}
\item[(1)] If $x \in X_{i} \setminus \{x_{i}\}$, then the neighborhood base
of $x$ in $\widetilde{\bigvee}_{i\in I}(X_i,x_i)$ is the one of $X_i$;
\item[(2)] The point $x^*$ has a neighborhood base, each element of which is
of the form:
$$\widetilde{\bigvee}_{i\in I \setminus F}(X_i,x_i) {\vee}
{\bigvee}_{j \in F}U_{j},$$
where $F$ is a finite subset of $I$ and each $U_j$
is an open neighborhood of $x_j$ in $X_j$ for $j \in F$.
\end{itemize}

\begin{lmm}\label{lmm:first}\cite[Theorem A.1]{E:free}
Suppose that the space
$X_i$ is
locally simply-connected and first countable at $x_i$,
for
 each $i\in I$. Then
 $$\pi _1(\widetilde{\bigvee}_{i\in I} (X_i,x_i),x^*)
\simeq \pmc{$\times$}\, \, \; _{i\in I} ^\si \pi _1(X_i,x_i).$$
\end{lmm}

\begin{lmm}\label{lmm:second}\cite[Proposition 2.10]{E:embed}
Let $X_i$ and $Y_j$ be locally simply-connected and first countable at
 $x_i$ and $y_j$, respectively for each $i\in I$ and $j\in J$. Then for the
 continuous map
 $$f: (\widetilde{\bigvee}_{i\in I}(X_i,x_i),x^*) \to
 (\widetilde{\bigvee}_{j\in J}(Y_j, y_j),y^*),$$
 the induced homomorphism
 $$f_*: \pi _1 (\widetilde{\bigvee}_{i\in I}(X_i,x_i),x^*) \to \pi _1
 (\widetilde{\bigvee}_{j\in J}(Y_j, y_j),y^*)$$
 is standard under the
 natural identifications:
\begin{eqnarray*}
\pi _1(\widetilde{\bigvee}_{i\in I} (X_i,x_i),x^*)
&=& \pmc{$\times$}\, \, \; _{i\in I} ^\si \pi _1(X_i,x_i) \ \
\hbox{and} \\
\pi _1(\widetilde{\bigvee}_{j\in J}(Y_j, y_j),y^*)
&=& \pmc{$\times$}\, \, \; _{j\in J} ^\si \pi _1(Y_j,y_j).
\end{eqnarray*}
\end{lmm}
Let $Y_0 = p^{-1}(\II\times [0,2/3))$ and $Y_1 =
p^{-1}(\II\times (1/3,1])$.
Then $SC(Z)= Y_0\cup Y_1$ and $Y_0\cap Y_1$ is open in $SC(Z)$.
We let $i_0: Y_0\cap Y_1 \to Y_0$, $i_1:Y_0\cap Y_1\to
 Y_1$, $j_0: Y_0\to SC(X)$, $j_1: Y_1\to SC(X)$, and $i:Y_0\cap Y_1\to
 SC(X)$ be the inclusion maps.
\medskip

\noindent
{\it Proof of\/} Theorem~\ref{thm:main1}.
We observe that $p^{-1}(\II\times \{ 1/2\})$ is a strong deformation
retract of $Y_0\cap Y_1$.
Let $C_n$ be the points on $\II \times \{ 1/2\}$
such that $C_{2n-1}$ is on the segment $A_nB_n$ and $C_{2n}$ is on the
segment $B_nA_{n+1}$.
Let $X_n$ be the subspace $[C,C_n]\cup p^{-1}(\{ C_n\})$ of $SC(Z)$.
Then $Y_0\cap Y_1$ is homotopy equivalent to
$\widetilde{\bigvee}_{n< \om}(X_n,C)$. Since $X_n$ is locally simply
connected and first countable at $C$ and $p^{-1}(\{ C_n\})$ is
homeomorphic to $Z$, $\pi _1(Y_0\cap Y_1)$ is isomorphic to
$\pmc{$\times$} \,\, \; _{n< \om} \pi _1(p^{-1}(\{ C_n\}) \cong  
\pmc{$\times$} \,\, \; _{n< \om} \pi _1(Z)$ 
by Lemmas~\ref{lmm:first} and
\ref{lmm:second}. Simlilarly, $\pi _1(Y_0)$ and $\pi _1(Y_1)$ are
isomorphic to $\pmc{$\times$}\, \, \; _{n< \om} \pi _1(p^{-1}(\{ A_n \})
\cong  
\pmc{$\times$} \,\, \; _{n< \om} \pi _1(Z)$
and $\pmc{$\times$}\, \, \; _{n< \om} \pi _1(p^{-1}(\{ B_n \})\cong  
\pmc{$\times$} \,\, \; _{n< \om} \pi _1(Z)$
respectively. 
Here we remark that $i_{0*}$ and $i_{1*}$ are standard homomorphisms
under these presentations of the fundamental groups. 

Since $Y_0, Y_1$ and $Y_0\cap Y_1$ are path-connected and open in
$SC(Z)$, we can apply the van Kampen theorem \cite[Theorem
2.1]{Massey:topology}
for homomorphisms $i_{0*}, i_{1*}, j_{0*}, j_{1*}$ and $i_*$ between
fundamental groups. The diagram formed by these five homomorphisms is a
pushout diagram and hence the ranges of $j_{0*}$ and $j_{1*}$ 
generate $\pi _1(SC(Z))$. Therefore $i_*$ is surjective. 
For the simple connectivity of $SC(Z)$ it suffices to show that $i_*$
is trivial. 

We let $u_n\in \pi _1(p^{-1}(\{ C_n\})$ be the copy of $u\in
\pi _1(Z)$.
Let $U_n$ be the word
$$u_n^{-1}u_{n+1}\cdots
u_{n +2k}^{-1}u_{n+ 2k+1}\cdots.$$
Since $i_{0*}(u_{2m}^{-1}u_{2m+1}) = e$ and $i_{1*}(u_{2m-1}^{-1}u_{2m})
= e$ and $i_{0*}$ and $i_{1*}$ are standard homomorphisms, 
$i_{0*}(U_{2m}) = e$, $i_{1*}(U_{2m}) = i_{1*}(u_{2m}^{-1})$,  
$i_{0*}(U_{2m-1}) = i_{0*}(u_{2m-1}^{-1})$, and $i_{1*}(U_{2m-1}) = e$. 

Now, an arbitrary element of $\pi _1(Y_0\cap Y_1)$ is expressed by a word
$$W\in \mathcal{W}(\pi _1(p^{-1}(\{ C_n\})): n\in \mathbb{N}).$$
For each letter $u_n\in \pi _1(p^{-1}(\{ C_n\}))$ for an odd $n$ appearing in
$W$, we insert $U_{n+1}$ successively to $u_n$ and form $W^*$. 
Since $$U_{n+1}\in\mathcal{W}(\pi_1(p^{-1}(\{ C_m\})): m\ge n ),$$
$W^*$ is actually
a word in $\mathcal{W}(\pi _1(p^{-1}(\{ C_n\})): n\in \mathbb{N})$.
We let $W_0$ be the word obtained by deleting all letters in $\bigcup
\pi _1(p^{-1}(\{ C_{2n-1}\}))$ from $W$. 
Since $i_{0*}$ and $i_{1*}$ are standard homomorphisms, $i_{0*}(W^*) =
i_{0*}(W)$ and $i_{1*}(W^*) = i_{1*}(W_0)$. 
Now $$W_0\in \mathcal{W}(\pi _1(p^{-1}(\{ C_{2n}\})): n\in \mathbb{N}).$$
We again insert $U_{n+1}$ for each letter $u_n$ appearing in
$W_0$ and form $W_0^*$. Then, by the symmetrical argument as above, we
conclude that $i_{1*}(W_0^*) = i_{1*}(W_0)$ and $i_{0*}(W_0^*) = e$. 
Now,
\begin{eqnarray*}
i_*(W)&=& j_{0*}\circ i_{0*} (W^*) = j_{1*}\circ i_{1*} (W_0) \\ 
      &=& j_{1*}\circ i_{1*} (W_0^*) = j_{0*}\circ i_{0*} (W_0^*) =
       j_{0*}(e) = e, 
\end{eqnarray*}
which imples that $\pi _1(SC(Z))$ is indeed
trivial.
\qed

\section{Proof of Theorem~\ref{thm:main2}}

For every group $G_i$,
$\Pi^\si _{i\in I}G_i$
is the subgroup of $\Pi
_{i\in I}G_i$ consisting of elements $u$ such that $\{ i\in I: u(i)\neq
0)\}$ is countable.
A space $X$ is called {\it semi-locally strongly contractible at}
$x\in X$, if
there exists an open neighborhood $U\subset X$ of $x$
such that there exists a contraction of $U$ in $X$ to $x$ which fixes $x$
(cf. \cite{EK:aloha}).

\begin{lmm}\label{lmm:aloha}\cite[Theorem 1.1]{EK:aloha}
Let $n\ge 2$ and let $X_i$ be a space which is 
$(n\! -\!1)$-connected semi-locally
strongly contractible at $x_i$ for each $i\in I$. Then
$$\pi _n(\widetilde{\bigvee}_{i\in I}(X_i,x_i),x^*)\cong
\Pi ^\si _{i\in I}\pi _n(X_i,x_i).$$
\end{lmm}

\noindent
{\it Proof of\/} Theorem~\ref{thm:main2}.
We shall use the Mayer-Vietoris sequence instead of the van Kampen theorem (as in the
preceding proof).
Consider the following Mayer-Vietoris homology exact sequence (over
$\mathbb Z$) for the triad $(SC(Z); Y_0, Y_1)$ from Section 2:
\[
H_n(Y_0\cap Y_1)\overset{i_{0*}+i_{1*}}\longrightarrow H_n(Y_0)\oplus
H_n(Y_1)\overset{ j_{0*}+ j_{1*}}\longrightarrow
H_n(SC(Z))\overset{\partial}\longrightarrow H_{n-1}(Y_0\cap Y_1).
\]
Since $Z$ is $(n\! -\! 1)$-connected, $Y_0\cap Y_1$ is also $(n\! -\!
1)$-connected, which implies
that
$H_{n-1}(Y_0\cap Y_1) = \{ 0\}$. Therefore
it suffices to show that $i_{0*}+i_{1*}$ is surjective.

Note
that $Y_0\cap Y_1$, $Y_0$ and $Y_1$ are simply connected and
that
$p^{-1}(\II\times \{1/2\})$, $p^{-1}(\II\times \{0\})$ and
$p^{-1}(\II\times \{ 1\})$ are strong deformation retracts of
$Y_0\cap Y_1$, $Y_0$ and $Y_1$,
respectively.
For the same reason as  explained in the first paragraph of the proof of Theorem~\ref{thm:main1},
the local properties required in Lemma~\ref{lmm:aloha}
for
 $Y_0\cap Y_1$, $Y_0$ and $Y_1$  are satisfied and we have:

\begin{eqnarray*}
H_n(Y_0\cap Y_1) &=& \pi _n(Y_0\cap Y_1) = \Pi _{m=1}^\infty
H_n(p^{-1}(\{ C_m\})), \\
H_n(Y_0) &=& \pi _n(Y_0) = \Pi _{m=1}^\infty H_n(p^{-1}(\{ A_m\})), \mbox{ and }\\
H_n(Y_1) &=& \pi _n(Y_1) = \Pi _{m=1}^\infty H_n(p^{-1}(\{ B_m\})),
\end{eqnarray*}
where $A_m, B_m, C_m$ are the points indicated in Figure 1.
Since $p^{-1}(\{ C_m\})$, $p^{-1}(\{ A_m\})$ and $p^{-1}(\{ B_m\})$ are
homeomorphic to $Z$, we can identify the homology groups of these spaces
with $H_n(Z)$.
Therefore
for $u\in \Pi _{m=1}^\infty H _n(p^{-1}(\{ C_m\}))$
\begin{eqnarray*}
i_{0*}(u)(1) &=& u(1), \\
i_{0*}(u)(m) &=& u(2m-1) + u(2m-2) \quad\mbox{ for }
m\ge 2, \\
i_{1*}(u)(m) &=& u(2m-1) + u(2m) \quad\mbox{ for }m\ge 1.
\end{eqnarray*}
For any given $v\in H_n(Y_0), w\in H_n(Y_1)$, define:
\begin{eqnarray*}
u(2m-1) &=& \Si _{k=1}^m v(k) - \Si _{k=1}^{m-1}w(k)
\mbox{ and } \\
u(2m) &=& \Si _{k=1}^m w(k) - \Si _{k=1}^{m}v(k).
\end{eqnarray*}
Then $i_{0*}(u) = v$ and  $i_{1*}(u) = w$ and hence $(i_{0*} +
i_{1*})(u) = v + w$.
We have thus
shown that $H_2(SC(Z),\mathbb{Z})$ is trivial and consequently $\pi _2(SC(Z))$ is also trivial by
the Hurewicz Theorem and Theorem~\ref{thm:main1}.
\qed
\begin{rmk}\label{rmk:difference}
{\rm
The proof of Theorem~\ref{thm:main2} in \cite{EKR:glasnik} was
along the same line as the
 proof of \cite[Theorem 1.1]{EKR:topsin}, which contains a procedure to
avoid $p^{-1}((0,1]\times \{ 1\})$.  The use of the Mayer-Vietoris
 sequence above makes it possible for us to skip this procedure,
 as does
 the use of the van Kampen theorem in the proof of Theorem~\ref{thm:main1}.
 When $Z$ is not simply-connected, we cannot avoid $p^{-1}((0,1]\times
 \{ 1\})$ for $\pi _2(SC(Z))$, which reflects the nontriviality of
 $\pi _2(SC(Z))$ in Theorem~\ref{thm:main2}.

We remark that the presentations of the homotopy groups,
 i.e. Lemmas~\ref{lmm:first}, \ref{lmm:second} and \ref{lmm:aloha}, are
 also useful to make proofs shorter. That is, our previous proofs
 implicitly contain the procedures used in the proofs of
 the Mayer-Vietoris sequence, the van Kampen theorem and the lemmas. }
\end{rmk}

\section{Acknowledgements}
The authors thank  the referee for several comments and suggestions. 
This research was supported by the
Slovenian Research Agency grants P1-0292-0101, J1-9643-0101 and
J1-2057-0101.
The first author was also
supported by the
Grant-in-Aid for Scientific research (C) of Japan
No. 20540097.

\providecommand{\bysame}{\leavevmode\hbox to3em{\hrulefill}\thinspace}

\end{document}